\documentclass[11pt, eqno, twoside]{article}

\usepackage{amsfonts} \usepackage{amsmath} \usepackage{amsthm}
\usepackage{amssymb} \usepackage{graphics} \usepackage{latexsym}
\usepackage{amscd} \usepackage{graphicx}\usepackage{amsfonts}

\theoremstyle{plain}
\newtheorem*{theorem}{Theorem}
\newtheorem*{corollary}{Corollary}
\newtheorem*{lemma}{Lemma}
\newtheorem*{proposition}{Proposition}

\theoremstyle{definition}
\newtheorem*{definition}{Definition}
\newtheorem*{example}{Example}%[section]
\newcommand{\al} {{\alpha}}  \newcommand{\be} {{\beta}}
  
\newcommand{\de} {{\delta}}  \newcommand{\De} {{\mit\Delta}}
\newcommand{\la} {{\lambda}}

\newcommand{\N} {{\mathbb N}}  \newcommand{\Z} {{\mathbb Z}}
\newcommand{\K} {{\mathbb K}}  
\newcommand{\bl} {{\mathbb L}}  %\newcommand{Fa} {{\mathfrak A}}
\newcommand{\lto}{\,\longrightarrow\,}

\newcommand{\md}{{\rm{\,mod\,}}}

\newcommand{\sci}{\scriptstyle}
\newcommand{\im}{{\rm Im\,}}

\begin{document}
\pagestyle{myheadings} \markboth{\centerline{\small{\sc On Weyl
Resolutions Associated to certain Frobenius Twists}}}
         {\centerline{\small{\sc Mihalis Maliakas}}}

\title{\bf On  Weyl resolutions associated to Frobenius twists}\footnotetext{

         \hspace{0.3cm}{\it Keywords and phrases:} Weyl filtration dimension, Weyl resolution,
         Schur algebra.

         \hspace{0.3cm}Research partially supported by ELKE, University of Athens.}

\author{Mihalis Maliakas}

\date{}

\maketitle

\begin{abstract}
We construct Weyl resolutions associated to certain Frobenius twists
of divided powers. Using these and other related complexes we obtain
the Weyl filtration dimension of the Schur algebras $S(2,r)$, a
result due to A.~Parker.
\end{abstract}

\section*{Introduction}

Let $\wedge_n$ be the subring of the polynomial ring $\Z{[x_1
,\dots, x_n ]}$ consisting of the symmetric polynomials. If $\la$ is
a partition, by $s_\la (x_1 ,\dots, x_n )$ we denote the associated
Schur function. The $s_\la (x_1 ,\dots, x_n )$, where $\la$ has at
most $n$ parts, form a $\Z$-basis of $\wedge_n$.

In \cite[pg 140]{M} the following identity in $\wedge_n$ is shown
\[
    h_k (x^2 _1 ,\dots, x^2 _n )=
    \begin{array}{c}\sum\limits^{k}_{j=0}\end{array}(-1)^j s_{(2k-j,j)}
    (x_1 ,\dots, x_n ),
\]
where $h_k (x_1 ,\dots, x_n )=\sum\limits_{1\le i_1 \le \cdots \le
i_k \le n} x_{i_1}\cdots x_{i_k}$ is the $k$-th complete symmetric
function.

In this paper we construct a complex of Weyl modules for the general
linear group $GL_n (K)$, where $K$ is an infinite field of
characteristic 2, whose homology is concetrated in degree 0 with the
following properties. The degree 0 homology is the first Frobenius
twist of the $k$-th divided power of the natural representation of
$GL_n (K)$. The terms of the complex are Weyl modules of $GL_n (K)$
that have formal characters $s_{(2k-j,j)}(x_1 ,\dots, x_n )$ and the
differentials are maps between Weyl modules in characteristic 2 so
that we have a representation theoretic realization of the above
identity (Thm 1.4).

The Weyl filtration dimension of the Schur algebra $S(n,r)$
corresponding to the homogeneous polynomial representations of $GL_n
(K)$ of degree $r$ is in general unknown. A.~Parker \cite{P1} has
determined the Weyl filtration dimension of $S(2,r)$ and $S(3,r)$.
Using the above and other related complexes we obtain an alternative
deduction of the Weyl filtration dimension of $S(2,r)$ (Section 2).

\section{Complexes associated to Frobenius twists}

\noindent\textbf{1.1} Throughout this paper, $K$ denotes a field of
characheristic $p>0$, $V$ a $K$-vector space with $\dim V=n<\infty$
and $G$ the general linear group $GL(V)$. We will use freely the
notation of \cite{A-B-W}. For example, the Weyl (respectively,
Schur) module of $G$ corresponding to a partition $\la$ will be
denoted by $K_\la V$ (respectively, $L_\la V$). The set of
partitions is denoted by $\wedge^+$. Let $\wedge^+ (n)=\{ \la \in
\wedge^+ |\la =(\la_1 ,\dots, \la_n )\}$ and if $r\ge 0$ let
$\wedge^+ (n,r)=\Big\{ \la \in \wedge^+ (n) |
\sum\limits^{n}_{i=1}\la_i =r\Big\}$. %The transpose of a partition
%$\la$ is denoted by $\tilde{\la}$.

We define here certain maps involving the Hopf algebra structure of
the exterior algebra $\wedge V$ of $V$ and the divided power algebra
$DV$ of $V$. Let $AV$ be one of $\wedge V$ or $DV$. We have a
natural grading $AV=\underset{a\ge 0}{\oplus}{} A_a V$, where for
all $a$, $A_a V=\overset{a}{\wedge}{}V$ or $A_a V=D_a V$. The
commultiplication of $AV$ will be denoted by $\De :AV \to AV \otimes
AV$ and the multiplication by $m:AV\otimes AV \to AV$. We will use
the same notation $\De$ to indicate the component $A_{a+b}V\to A_a
V\otimes A_b V$ of the restriction of $\De$ to $A_{a+b}V$ and
similarly for $m$.

If $a,b,t\ge 0$, define a map
\[
    \partial_t : A_a V\otimes A_b V \to A_{a+t}V\otimes A_{b-t}V
\]
as the composition $A_a V\otimes A_b V \xrightarrow{1\otimes
\De}{}A_a V\otimes A_t V\otimes A_{b-t}V\xrightarrow{m\otimes
1}{}A_{a+t}V\otimes A_{b-t}V$, where it is understood that if $b-t
<0$, then $A_{b-t}V=0$. Then $\partial_t$ is a map of $G$-modules
with respect to the diagonal action of $G$ on $A_a V\otimes A_b V$
and $A_{a+t}V\otimes A_{b-t}V$.

Suppose $a,b,s,t\ge 0$. It is easy to verify that the following
diagram commutes
\[
\begin{CD}
    A_a V\otimes A_b V @>\partial_s >> A_{a+s}V\otimes A_{b-s}V\\
    @V{\partial_t}VV @VV{\partial_t}V \\
    A_{a+t}V\otimes A_{b-t}V @>\partial_s >> A_{a+s+t}V\otimes
    A_{b-s-t}V
\end{CD}
\]
and moreover
\[
    \partial_s \circ \partial_t =\partial_t \circ \partial_s
    =\binom{s+t}{t}\partial_{s+t},
\]
where $\binom{s+t}{t}$ is the indicated binomial coefficient and
$\partial_{s+t}$ is the map $A_a V\otimes A_b V\to A_{a+s+t}V\otimes
A_{b-s-t}V$. Suppose in addition that $t$ satisfies $0<t<p$ and let
$\bar{t}=p-t$. Since the characteristic of $K$ is $p$, we have
\[
    \partial_{\bar{t}}\circ \partial_t =\partial_t \circ
    \partial_{\bar{t}}=0.\tag{1}
\]
\textbf{1.2} In this subsection, we will describe the maps which
will be the differentials of the complexes to be defined in
subsection 1.3. The differentials are maps between certain Weyl (or
Schur) modules in positive characteristic. Such maps have been
studied, among others, by Carter and Payne in \cite{C-P} where they
show that ${\rm Hom}_G (K_\la V,K_\mu V)\neq 0$ for certain pairs of
partitions $\la ,\mu$ where $\mu$ is obtained from $\la$ by
``raising'' some boxes of a row of $\la$ to a higher row. We will
need the special case where the partitions consist of two rows. In
order to make the Carter-Payne maps explicit, we give a different
proof of this special case. We remark that this particular case is
considerably simpler than the general situation.

\begin{lemma}%\label{le1.1}
Let $\la, \mu \in \wedge^+ (2)$ such that $\la =(a,b)$, $\mu
=(a+d,b-d)$ for some $d>0$. Suppose there is a positive integer $e$
such that $d<p^e$ and $a-b+d+1 \equiv 0\md p^e$. Then
\begin{itemize}
    \item[a)] the map $\partial_d : D_a V\otimes D_b V\to
    D_{a+d}V\otimes D_{b-d}V$ induces a map of $G$-modules $K_\la
    V\to K_\mu V$ which is nonzero provided $n\ge 2$,
    \item[b)] the map $\partial_d : \overset{a}{\wedge}{}V\otimes
    \overset{b}{\wedge}{}V\to \overset{a+b}{\wedge}{}V\otimes
    \overset{b-d}{\wedge}{}V$ induces a map of $G$-modules $L_\la
    V\to L_\mu V$ which is nonzero provided $a+d\ge n$.
\end{itemize}
\end{lemma}

\begin{proof}
a) We will employ the following notation. Under the comultiplication
map $D_{i+j}V\to D_i V\otimes D_j V$, the image of $x$ is denoted by
$\sum\limits_\al x_\al (i)\otimes x_\al (j)'$. Recall from
\cite[II]{A-B-W} the exact sequence
\[
    \bigoplus^{b-1}_{\ell =0} D_{a+b-\ell}V\otimes D_\ell
    V\xrightarrow{\;\Box\;}{}D_a V\otimes D_b
    V\xrightarrow{\;d'_\la\;}{}K_\la V\to 0,
\]
where the restriction of $\Box$ to $D_{a+b-\ell}V\otimes D_\ell V$
is the composition $D_{a+b-\ell}V\otimes D_\ell V\xrightarrow{\De
\otimes 1}{}D_a V\otimes D_{b-\ell}V\otimes D_\ell V
\xrightarrow{1\otimes m}{}D_\ell V\otimes D_bV$ and $d'_\la$ is the
map $d'_\la (V)$ of \cite[Def.\;II.1.3]{A-B-W}. Consider the
relation $D_{a+b-\ell}V\otimes D_\ell V$ of $K_\la V$, $0\le \ell
\le b-1$.

\medskip\noindent\textbf{Case 1.} Suppose $\ell \le b-d-1$.

\noindent We note that the $D_{a+b-\ell +t}V\otimes D_{\ell -t}V$,
where $t=0,\dots, \ell$, are relations for $K_\mu V$. For each
$t=0,\dots, d$, let $f^\ell _t :D_{a+b-\ell}V\otimes D_{\ell}V \to
D_{a+b-\ell +t}\otimes D_{\ell -t}V$ be the composite map
\[
    D_{a+b-\ell}V\otimes D_{\ell}V \xrightarrow{\partial_t}
    D_{a+b-\ell +t}V\otimes D_{\ell -t}V \xrightarrow{\binom{a-\ell
    +d}{d-t}}{}D_{a+b-\ell +t}V\otimes D_{\ell -t}V,
\]
where the map on the right is multiplication by the indicated
binomial coefficient. Put $f^\ell =\sum\limits^{\ell}_{t=0}f^\ell
_t$. We claim that the following diagram commutes
\[
\begin{CD}
    D_{a+b-\ell} V\otimes D_\ell V @>\Box >> D_{a}V\otimes D_{b}V\\
    @V{f^\ell}VV @VV{\partial_d}V \\
    \bigoplus\limits^{d}_{t=0}D_{a+b-\ell +t}V\otimes D_{\ell -t}V @>\Box >> D_{a+d}V\otimes
    D_{b-d}V.
\end{CD}
\]
Indeed, let $x\otimes y\in D_{a+b-\ell} V\otimes D_\ell V$. In one
direction we have
\begin{align*}
    \partial_d \Box (x\otimes y)=\begin{array}{l}\sum\limits^{d}_{i=0}\binom{a+i}{i}
    \sum\limits_{\al ,\be}\end{array}x_\al (a+i)y_\be (d-i)\otimes x_\al (b-\ell -i)'
    y_\be (\ell -d +i)',\tag{2}
\end{align*}
where the binomial coefficient comes from multiplication in $DV$. In
the other direction we have
\begin{align*}
    \Box f^\ell (x\otimes y)=
    &\begin{array}{l}\sum\limits^{d}_{t=0}\binom{a-\ell +d}{d-t}\sum\limits^{t}_{s=0}\binom{\ell
    -t+s}{s}\sum\limits_{\al ,\be}\end{array}x_\al (a+d-t+s)y_\be (t-s )\\
    &\hspace{2.6cm}\otimes x_\al (b-\ell -d+t-s)' y_\be (\ell -t+s)'.
\end{align*}
Let
\begin{align*}
    A_i =&\begin{array}{l}\sum\limits_{t-s=d-i} \binom{a-\ell +d}{d-t}\binom{\ell -t+s}{s}
    \sum\limits_{\al ,\be}\end{array}x_\al (a+d-t+s)y_\be (t-s)\\
    &\hspace{2.8cm}\otimes x_\al (b-\ell -d +t-s)' y_\be (\ell -t +s)'
\end{align*}

\noindent so that $\Box f^\ell (x\otimes
y)=\sum\limits^{d}_{i=0}A_i$. We have
\begin{align*}
    A_i =
    \begin{array}{l}\sum\limits^{i}_{s=0} \binom{a-\ell +d}{i-s}\binom{\ell -d+i}{s}
    \sum\limits_{\al ,\be}\end{array}x_\al (a+i)y_\be (d-i)
    \otimes x_\al (b-\ell
    -i)' y_\be (\ell -d +i)'.\tag{3}
\end{align*}
Using the binomial identity $\sum\limits^{i}_{s=0}\binom{a-\ell
+d}{i-s}\binom{\ell -d+i}{s}=\binom{a+i}{i}$, it follows from (2)
and (3) that $\sum\limits^{d}_{i=0}A_i =\partial_d \Box (x\otimes
y)$ and hence $\Box f^\ell (x\otimes y)=\partial_d \Box (x\otimes
y)$. (In this case we did not use the assumption on $p$).

\medskip\noindent\textbf{Case 2.} Suppose $\ell >b-d+1$.

\noindent For each $t$ satisfying $\ell -b+d+1 \le t\le d$, let
$g^\ell _t$ be the composition
\[
    D_{a+b-\ell}V\otimes D_{\ell}V \xrightarrow{\partial_t}{}
    D_{a+b-\ell +t}V\otimes D_{\ell -t}V
    \xrightarrow{\binom{b-\ell -1}{d-t}}{}
    D_{a+b-\ell +t}V\otimes D_{\ell -t}V
\]
and put $g^\ell =\sum\limits^{d}_{t=c}g^\ell _t$ where $c=\ell
-b+d+1$. We claim that the following diagram commutes
\[
\begin{CD}
    D_{a+b-\ell} V\otimes D_\ell V @>\Box >> D_{a}V\otimes D_{b}V\\
    @V{g^\ell}VV @VV{\partial_d}V \\
    \bigoplus\limits^{d}_{t=c}D_{a+b-\ell +t}V\otimes D_{\ell -t}V @>\Box >> D_{a+d}V\otimes
    D_{b-d}V.
\end{CD}
\]
Indeed, let $x\otimes y\in D_{a+b-\ell}V\otimes D_{\ell}V$. In one
direction, we have equation (2). For the other direction we compute
\begin{align*}
    \Box g^\ell (x\otimes y)=&
    \begin{array}{l}\sum\limits^{d}_{t=c} \binom{b-\ell -1}{d-t}
    \sum\limits^{t}_{s=0}\binom{\ell -t+s}{s}
    \sum\limits_{\al ,\be}\end{array}
    x_\al (a+d -t+s) y_\be (t-s)\\
    &\hspace{2.6cm}\otimes x_\al (b-\ell -d+t-s)' y_\be
    (\ell -t+s)' .\tag{4}
\end{align*}
Note that if $t<c$, then $d-t>b-\ell -1$ so that $\binom{b-\ell
-1}{d-t}=0$. Thus in the right hand side of (4), the summation with
respect to $t$ may be taken from $t=0$ to $t=d$. Hence $\Box g^\ell
(x\otimes y)=\sum\limits^{d}_{i=0}B_i$, where
\begin{align*}
 B_i& =
    \begin{array}{l}
    \sum\limits_{t-s=d-i}\binom{b-\ell -1}{d-t}\binom{\ell
    -t+s}{s} \sum\limits_{\al ,\be}\end{array} x_\al (a+d-t+s)y_\be
    (t-s)\\
    &\hspace{3.4cm}\otimes x_\be (b-\ell -d+t-s)' y_\be (\ell -t+s)'\\
    &=\begin{array}{l}\sum\limits^{i}_{s=0}\binom{b-\ell -1}{i-s}\binom{\ell
    -d+i}{s}\sum\limits_{\al ,\be}\end{array}x_\al (a+i)y_\be (d-i)
    \otimes x_\al (b-\ell -i)' y_\be (\ell -d+i)'\\
    &=\begin{array}{l}\binom{b-1-d+i}{i}\sum\limits_{\al ,\be}\end{array}x_\al (a+i)y_\be
    (d-i)\otimes x_\al (b-\ell -i)' y_\be (\ell -d +i)'\tag{5}
\end{align*}
by our familiar binomial identity. From the assumptions on $d$ and
$p$ in the statement of the Lemma it follows that for each
$i=0,\dots,d$
\[
    \begin{array}{c}\binom{a+i}{i}\equiv \binom{b-1-d+i}{i}\md p.
    \end{array}
\]
From this and equations (2), (5) we see that the diagram commutes.

We have shown that the map $\partial_d : D_a V\otimes D_d V\to
D_{a+d}V\otimes D_{b-d}V$ induces a map $K_\la V\to K_\mu V$. If
$n\ge 2$, this last map is nonzero since the image of the standard
tableau $1^{(a)}|2^{(b)}$ is the tableau $1^{(a)}2^{(d)}|2^{(b-d)}$
which is standard because $a\ge b-d$.

b) The proof is similar to the proof of a) and thus omitted.\\
\end{proof}

We will denoted again by $\partial_d$ the maps $K_\la V\to K_\mu V$
and $L_\la V\to L_\mu V$ provided by the previous Lemma.

\bigskip\noindent\textbf{1.3} In this section, we describe a family
of complexes of $G$-modules associated to certain Frobenius twists.

We adopt the convention that $K_\al V=L_\al V=0$ if $\al$ is a
sequence of non negative integers that is not a partition. This is
not used in \cite{A-B-W} but is in line with the representation
theory of algebraic groups \cite[II 2.6 Prop.]{J}.

\begin{definition}%\label{de1.2}
Suppose \ $d,r$ \ are \ positive \ integer \ such  \ that \ $0<d<p$
\ and \ $r-d+1 \equiv 0\md p$. Let $\K_* (r,d,V)$ be the complex
where for each $i\ge 0$
\[
    \K_{2i+1}(r,d,V)=K_{(r-ip-d,ip+d)}V, \quad
    \K_{2i}(r,d,V)=K_{(r-ip,ip)}V
\]
and the differential $\de_i : \K_i (r,d,V)\to \K_{i-1}(r,d,V)$ is
given by
\[
    \de_{2i+1}=\partial_d , \quad \de_{2i}=\partial_{\bar{d}},
\]
where $\bar{d}=p-d$.
\end{definition}

We remark that the maps $\de_i$ in the above definition are provided
by the conclusion of Lemma 1.2 and they alternate between
``raising'' $d$ and $p-d$ boxes. Moreover, $\K_* (r,d,V)$ is a
complex because of (1). If $j>\lfloor \frac{r}{p}\rfloor$, where
$\lfloor\frac{r}{p}\rfloor$ is the integer part of $\frac{r}{p}$,
then according to our convention mentioned above, we have $\K_j
(r,d,V)=0$.

In a similar manner we define complexes $\bl_* (r,d,V)$ by replacing
the Weyl modules $K_{(r-ip-d,ip+d)}V$, $K_{(r-ip,ip)}V$ by the Schur
modules $L_{(r-ip-d,ip+d)}V$, $L_{(r-ip,ip)}V$ respectively.

\begin{example}%\label{ex1.3}
If $p=2$ (and hence $d=1$ and $r$ is even), then $\K_* (r,1,V)$ is
the complex
\[
    0\lto K_{(r/2,r/2)}V\xrightarrow{\,\partial_1\,}{}\cdots
    \lto K_{(r-2,2)}V\xrightarrow{\,\partial_1\,}{}K_{(r-1,1)}
    V\xrightarrow{\,\partial_1\,}{}D_r V
\]
and $\bl_* (r,1,V)$ is
\[
    0\lto L_{(r/2,r/2)}V\xrightarrow{\,\partial_1\,}{}\cdots
    \lto L_{(r-2,2)}V\xrightarrow{\,\partial_1\,}{}L_{(r-1,1)}
    V\xrightarrow{\,\partial_1\,}{}\overset{r}{\wedge} V.
\]
\end{example}

\noindent\textbf{1.4} Let $n=\dim V$. In this section we determine
the homology of the complexes $\K_* (r,d,V)$, $\bl_* (r,d,V)$ in the
cases 1) $p=2$ (and any $n$) and 2) $n=2$ (and any $p$). The $m$-th
Frobenius twist of a $G$-module $M$ will be denoted $M^{(m)}$.

\begin{theorem}%\label{th1.4}
Suppose $d,r$ are positive integers such that $0<d<p$ and
$r-d+1\equiv 0\md p$. If $p=2$ or $n=2$, then $\K_* (r,d,V)$ has
homology concentrated in degree $0$ and $\bl_* (r,d,V)$ has homology
concentrated in degree $\lfloor\frac{r}{p}\rfloor$. In particular,
for $p=2$ (and hence $d=1$ and $r$ is even, $r=2r'$) we have
\[
    H_0 (\K_* (r,1,V))=(D_{r'}V)^{(1)}, \quad H_{r'}(\bl_*
    (r,1,V))=(\overset{r'}{\wedge}V)^{(1)}.
\]
\end{theorem}

\begin{proof}
a) Suppose $p=2$.\\
For the bookkeeping to follow, we observe that if $m\ge \ell$ and
$W$ is a $K$-vector space of dimension 1, then the complex
\begin{align*}
   & C_* (m,\ell ):\\
   &0\lto D_{m-\ell}W\otimes D_\ell W\xrightarrow{\partial_1}
    D_{m-\ell +1}W\otimes D_{\ell -1}W\xrightarrow{\partial_1}\cdots\;\lto
    D_{m-1}W\otimes D_1 W\xrightarrow{\partial_1}D_{m}W
\end{align*}
is isomorphic to
\[
    0\lto K\xrightarrow{m-\ell +1}K\xrightarrow{m-\ell +2}K\lto
    \cdots \xrightarrow{\,m\,}K
\]
and thus: 1) if both $m, \ell$ are odd, then $C_* (m,\ell)$ is
exact, 2) if both $m,\ell$ are even, then $C_* (m,\ell)$ has
homology concentrated in degree 0 and this is $D_m W$.

We prove by induction on $n$ that the homology of $\K_* (r,1,V)$ is
concentrated in degree 0.

The case $n=1$ being trivial, let $n\ge 2$ and suppose $V=U\oplus W$
where $U,W$ are subspaces of $V$ with $\dim W=1$. Fix a basis $x_1
,\dots, x_{n-1}$ of $U$ and $x_n$ of $W$ and consider he ordering
$x_1 <x_2 <\cdots < x_{n-1}<x_n$. By the analog of
\cite[Thm.\;II.4.11]{A-B-W} for Weyl modules, for each $i$
satisfying $2i\le r$ we have the following direct sum decomposition
of vector spaces
\[
    K_{(r-i,i)}V\simeq \underset{s,t}{\oplus}K_{(k-i-s,i-t)}U\otimes
    D_s W \otimes D_t W,
\]
where $0\le s\le r-2i$, $0\le t\le i$. Under this identification the
differential of $\K_* (r,1,V)$ looks like

\begin{minipage}{12cm}
\[
\begin{picture}(400,60)
 \put(10,50){\makebox{$K_{(r-i-s,i-t)} U\otimes D_s W\otimes D_t W$}}
 \put(158,55){\makebox{$\scriptstyle{\al \otimes 1\otimes 1}$}}
 \put(155,53){\vector(1,0){35}}
 \put(193,50){\makebox{$K_{(r-i-s+1,i-t-1)}U\otimes D_s W\otimes D_t
 W$}}
 \put(193,13){\makebox{$K_{(r-i-s,i-t)}U\otimes D_{s+1} W\otimes
 D_{t-1}W$}}
 \put(155,50){\vector(1,-1){33}}
 \put(155,25){\makebox{$\scriptstyle{1 \otimes \be}$}}
 \end{picture} %\tag{6}
\]
\end{minipage}\hfill(6)

\noindent where $\al : K_{(r-i-s,i-t)}U\lto K_{(r-i-s+1,i-t-1)}U$
sends a standard tableau

%$1^{(a_1 )}\cdots (n-1)^{(a_n )}| 1^{(b_1 )}\cdots n^{(b_n )}$ to
%
%$\sum\limits^{n}_{j=1}(a_1 +1)1^{(a_1 )}\cdots j^{(a_j +1)}\cdots
%n^{(a_n )}|1^{(b_1 )}\cdots j^{(b_j -1)}n^{(b_n )}$

\noindent$1^{(a_1 )}\cdots (n-1)^{(a_{n-1} )}| 1^{(b_1 )}\cdots
(n-1)^{(b_{n-1})}$ to

\noindent$\sum\limits^{n-1}_{j=1}(a_j +1)1^{(a_1 )}\cdots j^{(a_j
+1)}\cdots (n-1)^{(a_{n-1})}|1^{(b_1 )}\cdots j^{(b_j
-1)}(n-1)^{(b_{n-1})}$

\noindent and $\be : D_s W\otimes D_t W\to D_{s+1}W\otimes D_{t-1}W$
is the differential of $C_*$.

Since the differential of $\K_* (r,1,V)$ is a map of $G$-modules, it
is a map of $T$-modules where $T$ is the maximal torus of diagonal
matrices in $GL_n (K)\simeq GL(V)$. Thus we may consider the
subcomplex $\K_* (r,V;m)$ of $\K_* (r,1,V)$ of $x_n$-content of
equal to $m$, that is, $\K_i (r,V;m)$ is the subspace of
$K_{(r-i,i)}V$ with basis the standard tableaux of $K_{(r-i,i)}V$ in
which the element $x_n$ appears $m$ times. We have
\[
    \K_* (r,1,V)=\underset{m\ge 0}{\oplus}\K_* (r,V;m).
\]
Also
\[
    \K_i (r,V;m)=\underset{s,t}{\oplus} K_{(r-i-s,i-t)}U\otimes D_s
    W\otimes D_t W\tag{7}
\]
where the sum ranges over all $s,t$ subject to the conditions $0\le
s\le r-2i$, $0\le t\le i$, $s+t=m$.

\medskip\noindent\textbf{Claim 1.} Let $m$ be odd. Then $\K_*
(r,V;m)$ is exact.

Indeed, since $r-m$ is odd, let $r-m=2q+1$. From (6) and (7) it
follows that $\K_* (r,V;m)$ is the total complex of the following
double complex.
\[
\begin{array}{crccc}
\sci{K_{(q+1,q)}U\otimes D_m W }&\to \cdots
\to&\sci{K_{(r-m-1,1)}U\otimes D_m
 W}&\to&\sci{K_{(r-m)}U\otimes D_m W} \\
 \sci\uparrow & & \uparrow & & \uparrow \\
 \sci{K_{(q+1,q)}U\otimes D_{m-1} W \otimes D_1 W} &\to \cdots \to&\sci{K_{(r-m-1,1)}U\otimes
 D_{m-1}\otimes D_1 W }&\to&\sci{K_{(r-m)}U\otimes D_{m-1} W\otimes D_1 W} \\
 \uparrow & & \uparrow & & \uparrow \\
 \sci{0} & \cdots \to & \sci{K_{(r-m-1,1)}U\otimes
 D_{m-2}\otimes D_2 W }&\to&\sci{K_{(r-m)}U\otimes D_{m-2} W\otimes D_2 W} \\
  & & \uparrow  & & \uparrow \\
  & & \vdots & & \vdots \\
  & &  & & \uparrow \\
 & \cdots \to & {\cdots}%{K_{(r-m-1,1)}U\otimes D_{m-\ell_2}\otimes D_{\ell_2} W }
 &\to&\sci{K_{(r-m)}U\otimes D_{m-\ell_0} W\otimes D_{\ell_0} W} \\
  & & \uparrow & & \uparrow \\
  & & \sci{0} & & \sci{0}
 \end{array}
\]
The $i$-th column, $i=0,1,\dots$ (counting from right to left) for
$i$ satisfying $r-m-2i \ge 0$ is the complex $K_{(r-m-i,i)}U\otimes
C_* (m,\ell_i )$, where $\ell_i =\min \{ m,r-m-2i\}$. Each nonzero
horizontal map is of the form $\al \otimes 1\otimes 1$ (see (6)).
Now since both $m$ and $\ell_i$ (and $\ell_i$) are odd, it follows
from the observation made at the beginning of the proof that each
column of the double complex is exact. Hence $\K_* (r,V;m)$ is
exact.

\medskip\noindent\textbf{Claim 2.} Let $m$ be even. Then $\K_*
(r,V;m)$ has homology concentrated in degree zero.

Indeed, since $r-m$ is even, let $r-m=2q$. Again $\K_* (r,V;m)$ is
the total complex of the double complex whose $i$-th column is
$K_{(r-m-i,i)}U\otimes C_* (m,\ell_i )$ where $\ell_i =\min \{
m,r-m-2i\}$. Notice however, with respect to the display of the
bicomplex in the previous case, that the left most column now has
the form
\[
    \begin{array}{c}
    \sci{K_{(q,q)}U\otimes D_m W}\\
    \uparrow\\
    \sci{0}
    \end{array}
\]
Since both $m$ and $\ell_i$ are even, we know from the observation
made earlier that for every $i$, the $i$-th column has homology
concentrated in degree 0 and this is $K_{(r-m-i,i)}U\otimes D_m W$.
Thus the $E_1$ page of the first spectral sequence associated to our
bicomplex is the top row
\[
    0\lto K_{(q,q)}U\otimes D_m W\lto \cdots \lto K_{(r-m-1,1)}
    U\otimes D_m W\lto K_{(r-m)}U\otimes D_m W
\]
of the bicomplex, which is $\K_* (r-m,1,U)\otimes D_m W$. By the
inductive hypothesis, this complex has homology concentrated in
degree 0. Thus the spectral sequence collapses on both axes and
hence $\K_* (r,V;m)$ has homology concentrated in degree zero.

From Claim 1 and Claim 2 it follows that $\K_* (r,1,V)$ has homology
concentrated in degree zero.

We consider now $H_0 (\K_* (r,1,V))$. Let $x_1 ,\dots, x_n$ be a
basis of $V$ and consider a tableau $x=x^{(a_1)}_1 \cdots
x^{(a_n)}_n |x_i \in K_{(r-1,1)}V$. For the map $\partial_1
:K_{(r-1,1)}V\to D_r V$ we have $\partial_1 (x)=(a_i +1)x^{(a_1)}_1
\cdots x^{(a_i +1)}_i \cdots x^{(a_n)}_n$. Since the characteristic
of $K$ is 2, it follows that a basis of $H_0 (\K_* (r,1,V)$ is
\[
    \{ x^{(2a_1)}_1 \cdots x^{(2a_n)}_n +{\rm Im} \partial_1 | a_1
    +\cdots + a_n =r' \}
\]
Now it is straightforward to verify that the vector space
isomorphism
\[
    (D_{r'}V)^{(1)} \ni x^{(a_1)}_1 \cdots x^{(a_n)}_n
    \mapsto x^{(2a_1)}_1 \cdots x^{(2a_n)}_n +{\rm Im}\partial_1 \in
    H_0 (\K_* (r,1,V))
\]
is $GL_n (K)$-equivariant.

Next we consider the complex $\bl_*(r,1,V)$. Most of the details of
the proof are similar (and in fact simpler) to the previous proof
and for these we will be consice. We argue by induction on $n=\dim
V$ that $\bl_* (r,1,V)$ has homology concentrated in degree $r'$,
the case $n=1$ being trivial. Let $n\ge 2$ and, as before, let
$V=U\oplus W$, $\dim W=1$. Let $x_1 ,\dots, x_{n-1}$ be a basis of
$U$ and $x_n$ a basis of $W$. From \cite[Thm.\;II.4.11]{A-B-W} we
have $L_{(r-i,i)} V=L_{(r-i,i)} U\oplus L_{(r-i-1,i)} U\otimes
W\oplus L_{(r-i,i-1)} U\otimes W\oplus L_{(r-i-1,i-1)}U\otimes
W\otimes W$.

Let $\bl_* (r,V;m)$ be the subcomplex of $L_* (r,1,V)$ of
$x_n$-content equal to $m$. From the standard basis theorem for
Schur modules \cite[Thm.\;II.2.16]{A-B-W} it follows that $\bl_*
(r,V;m)=0$ for $m\ge 3$. Hence $\bl_* (r,1,V)
=\bigoplus\limits^{2}_{m=0}\bl_* (r,V;m)$.

For $m=1$, we see that $\bl_* (r,V;1)$ is the total complex of the
double complex
\[
\begin{array}{crccc}
{L_{(r',{r'}-1)}U\otimes W}&\to \cdots \to&{L_{(r-2,1)}U\otimes
W }&\to&{L_{(r-1)}U\otimes W} \\
\uparrow & & \uparrow & & \uparrow \\
{L_{(r',{r'}-1)}U\otimes W}&\to \cdots \to&{L_{(r-2,1)}U\otimes
W }&\to&{L_{(r-1)}U\otimes W}. \\
 \end{array}
\]
Since the vertical maps are identity maps, the total complex is
exact.

For $m=0,2$ we see that
\begin{align*}
    &\bl_* (r,V;0) =\bl_* (r,1,U),\\
    &\bl_* (r,V;2)=\bl_* (r-2,1,U)[-1]\otimes W\otimes W
\end{align*}
where $\bl_* (r-2,1,U)[-1]$ is $\bl_* (r-2,1,U)$ with degrees
shifted by $-1$. By induction, both of the above complexes have
homology concentrated in degree $r'$.

Next we consider $H_{r'}(\bl_* (r,1,V))=\ker (L_{(r',
r')}\xrightarrow{\partial_1}{}L_{(r' +1,r' -1)}V)$ which we denote
by $M$. Let $x_1 ,\cdots, x_n$ be a basis of $V$ and assume $n\ge
r'$ so that $L_{(r',r')}V\neq 0$. Let $x\in L_{(r' ,r')}V$ be the
standard tableau $x=x_1 \cdots x_{r'} | x_1 \cdots x_{r'}$. Then
$x\in M$. Let us fix the maximal torus of diagonal matrices in $GL_n
(K)$ and the Borel subgroup of lower triangular matrices in $GL_n
(K)$. Then for any partition $\la$, $L_\la V$ is the induced module
of $GL_n (K)$ of highest weight $\tilde{\la}$ \cite{D}, where
$\tilde{\la}$ denotes the transpose partition of $\la$, and its
formal character is the Schur function $s_{\tilde{\la}}(x_1 ,\dots,
x_n )$. Now $x\in M$ is a highest weight vector of $L_{(r' ,r')}V$
of weight $(2^{r'})=(2,\dots, 2)$. Since $x\in M$ and $L_{(r'
,r')}V$ has simple socle generated by its highest weight vector
\cite{J}, we conclude that $L(2^{r'})\subseteq M$, where $L(2^{r'})$
is the irreducible $GL_n (K)$-module of highest weight $(2^{r'})$.
We have $L(2^{r'})=L(1^{r'})^{(1)}=(\overset{r'}{\wedge}V)^{(1)}$,
so that
\[
    (\overset{r'}{\wedge}V)^{(1)}\subseteq M.\tag{7}
\]
In order to show equality in (7) we will prove that $\dim
M=\binom{n}{r'}$.

From the complex $\bl_* (r,1,V)$ that has homology $M$ concentrated
on the extreme left we have
\[
    \dim M=\begin{array}{l}\sum\limits^{r'}_{j=0}\end{array}
    \dim L_{(r' +j,r' -j)}V.
\]
It is well know that $\dim L_{(r' +j,r' -j)}V= \binom{n}{r' +j}
\binom{n}{r' -j}-\binom{n}{r' +j+1}\binom{n}{r' -j-1}$
\cite[I\;(3.5)]{M}. By substituting we have
\[
    \begin{array}{l}
    \dim M=2\sum\limits^{r'}_{j=0}(-1)^j \binom{n}{r' +j}
    \binom{n}{r' -j} -\binom{n}{r'}^2\end{array}.\tag{8}
\]
Using the identity
\[
    \begin{array}{l}
    2\sum\limits^{r'}_{j=0} (-1)^j \binom{n}{r' +j} \binom{n}{r' -j}
    =\binom{n}{r'}+\binom{n}{r'}^2 \end{array},
\]
which may be proven easily by induction on $n$, (8) yields $\dim
M=\binom{n}{r'}$.

b) Suppose $n=2$.\\
Let $x,y$ be a basis of $V$. From the definition of $\K_* (r,d,V)$,
it follows that $\ker \de_{2i+1}$ is generated by the standard
tableaux $x^{(a)}y^{(b)}|y^{(c)}\in K_{(r-2ip-d,2ip+d)}V$ such that
\[
    \begin{array}{l}\binom{b+d}{d}\equiv 0\md p\end{array},
\]
while the image of $\de_{2i+2}$ is generated by the standard
tableaux $x^{(a)}y^{(b)}|y^{(c)}$ such that
\[
    \begin{array}{l}\binom{b}{p-d}\not\equiv 0\md p\end{array}.
\]
Since among $p$ consecutive integers exactly one is divisible by
$p$, we obtain $\im \de_{2i+2}=\ker \de_{2i+1}$ for all $i\ge 0$. In
a similar way we have $\im \de_{2i+1}=\ker \de_{2i}$ for all $i\ge
1$, so that $\K_* (r,d,V)$ has homology concentrated in degree 0.\\
\end{proof}

\medskip%\noindent\textbf{1.5} Suppose $n=2$. Then $\K_* (r,d,V)$ has
%homology concentrated in degree 0. For $d=1$, whence $r$ is
%divisible by $p$, this homology is the cokernel of the map
%$K_{(r-1,1)}V \xrightarrow{\partial_1}{}D_r V$ and it follows easily
%that
%\[
%    H_0 (\K_* (r,1,V))=(D_{r'}V)^{(1)}.
%\]
\noindent\textbf{1.5} We note here two complexes related to those of
(1.3) that will be needed later. Suppose $r$ is odd and $p=2$.

1) Let $r\equiv 1\md 4$. By Lemma {1.2}, we have a nonzero map of
Weyl modules $K_{(r-2i,2i)} V\lto K_{(r-2i+2,2i-2)}V$, $1\le i\le
\frac{r}{4}$, induced by the map $D_{r-2i}V\otimes
D_{2i}V\xrightarrow{\partial_2}{}D_{r-2i+2}V\otimes D_{2i-2}V$. The
complex
\[
    \cdots \xrightarrow{\partial_2} K_{(r-4,4)}V
    \xrightarrow{\partial_2} K_{(r-2,2)}V\xrightarrow{\partial_2} D_r V
\]
will be denoted $\mathbb{M}_{*} (r,2,V)$. Note that its length is
$\lfloor \frac{r}{4}\rfloor$. By a direct computation (as in the
last part of 1.4), we see that if $n=2$, then $\mathbb{M}_{*}
(r,2,V)$ has homology concentrated in degree 0.

2) Let $r\equiv 3\md 4$. By Lemma 1.2 we have a nonzero map of Weyl
modules $K_{(r-2i-1,2i+1)}V\to K_{(r-2i+1,2i-1)}V$, $1\le i\le
\lfloor \frac{r}{4}\rfloor$, induced by the map $D_{r-2i-1}V\otimes
D_{2i+1}V\xrightarrow{\partial_2}D_{r-2i+1}V\oplus D_{2i-1}V$. The
complex
\[
    \cdots \xrightarrow{\partial_2} K_{(r-5,5)} V
    \xrightarrow{\partial_2} K_{(r-3,3)}
    \xrightarrow{\partial_2} K_{(r,1)}V
\]
will be denoted $\mathbb{N}(r,2,V)$. Note that its length is
$\lfloor \frac{r}{4}\rfloor$. By a direct computation, we see that
if $n=2$, then $\mathbb{N}_{*}(r,2,V)$ has homology concentrated in
degree 0.

\section{The Weyl filtration dimension of $S(2,r)$}

The Weyl filtration dimension of the Schur algebra $S(2,r)$ of $GL_2
(K)$ has been determined by Parker in \cite{P1}. See also \cite{P2}
for a different proof. In this section we utilise the resolutions of
Theorem 1.4 and (1.5) to obtain an alternative proof of Parker's
result.

\medskip\noindent\textbf{2.1} We consider the Schur algebra $S(n,r)$
corresponding to the homogeneous polynomial representations of $GL_n
(K)$ of degree $r$ \cite{G}. By $\md S(n,r)$ we denote the category
of finite dimensional $S(n,r)$-modules. Let $V$ be the natural $GL_n
(K)$-module. For each $\la \in \wedge^+ (n,r)$, we have the induced
module $\nabla(\la)$, the Weyl module $\triangle(\la)$ and the
irreducible module $L(\la )$, all of highest weight $\la$. It is
well known that $\nabla(\la) \simeq L_{\tilde{\la}}V$ and
$\triangle(\la)\simeq K_\la V$ \cite{D}.

We recall here relevant definitions. Let $X\in \md S(n,r)$.

a) We say that $X$ admits a Weyl filtration if there is a filtration
$0=X_d \subseteq \cdots \subseteq X_0 =X$ of $X$ by
$S(n,r)$-submodules such that each $X_i /X_{i+1}$ is either 0 or
isomorphic to some $K_\la V$, $\la \in \wedge^+ (n,r)$.

b) A resolution of $S(n,r)$-modules of the form
\[
    0\lto X_d \lto \cdots \lto X_0 \lto X\lto 0,
\]
where each $X_i$ admits a Weyl filtration, is called a Weyl
resolution of $X$ of length $d$.

c) If $X$ has a Weyl resolution of length $d$ and has no Weyl
resolution of length less than $d$, we say that the Weyl filtration
dimension of $X$ is $d$ (abbreviated $wfd(X)=d$).

d) The supremum of $\{ wfd(X)|X\in \md S(n,r)\}$ is called the Weyl
filtration dimension of $S(n,r)$ and is denoted by $wfd(S(n,r))$.
(This should not be confused with the Weyl filtration dimension of
$S(n,r)$ as a module over itself).

In the next proposition, we collect the well known facts concerning
Weyl filtration dimensions that will be needed. We write ${\rm
Ext}^i (-,-)$ in place of ${\rm Ext}^i _{GL_n (K)}(-,-)$.

\begin{proposition}%\label{pro2.1} $ $
\begin{itemize}
    \item[i)] Let $\la ,\tilde{\mu}\in \wedge^+ (n,r)$. Then
    \[
    {\rm Ext}^i (K_\la V,L_{\tilde{\mu}}V)=\left\{\begin{array}{ll}
    K,& \text{if } \ i=0 \ \text{ and } \ \la ={\mu}\\
    0, &\text{otherwise.}
    \end{array}\right.
    \]
    \item[ii)] Let $d$ be a non negative integer and $X\in \md
    S(n,r)$. Then the following are equivalent
    \begin{itemize}
        \item[a)] $wfd(X)=d$
        \item[b)] ${\rm Ext}^i (X,L_{\tilde{\la}}V)=0$ for all $i>d$
        and all $\la \in \wedge^+ (n,r)$ and ${\rm Ext}^d
        (X,L_{\tilde{\mu}}V)\neq 0$ for some $\mu \in \wedge^+
        (n,r)$.
    \end{itemize}
    \item[iii)] If $X\in \md S(n,r_1 )$, $Y\in \md  S(n,r_2 )$ have
    Weyl filtrations, then so does
    $X\otimes Y$.
    \item[iv)] If $X\in \md S(n,r_1 )$, $Y\in\md S(n,r_2 )$, then $wfd(X\otimes Y)\le wfd (X)+wfd (Y)$.
    \item[v)] $wfd (S(n,r))=\max \{ wfd (L(\la))| \la \in \wedge^+
    (n,r)\}.$
    \item[vi)] If \;$0\to X_d \to \cdots \to X_0 \to X\to 0$
is an exact sequence in $\md S(n,r)$, then $wfd(X)\le \max \{ wfd
(X_i )+i |i=0,\dots, d\}$.
\end{itemize}
\end{proposition}

For the proofs of i), ii) and iv) see \cite{F-P} (where the dual
notion of good filtration dimension was considered) and for iii) see
\cite{W}. The proofs of v) and vi) are standard using ii) and the
long exact sequence in cohomology.

\medskip\noindent\textbf{2.2} Suppose $Y\in \md S(n,r)$ has a Weyl
filtration. From part i) of the above proposition, it follows that
the multiplicity of $K_\la V$ as a factor of $Y$ is $\dim {\rm
Hom}(Y,L_{\tilde{\la}}V)$ (and hence does not depend on the choice
of good filtration).

\begin{lemma}%\label{le2.2}
Let $X\in \md S(n,r)$. If
\[
    0\lto X_d \lto X_{d-1} \lto \cdots \lto X_0 \lto X \lto 0,
\]
is a Weyl resolution of $X$ such that $X_d$ has a factor which is
not a factor of $X_{d-1}$, then $wfd(X)=d$.
\end{lemma}

\begin{proof}
Clearly we have $wfd(X)\le d$. Let $Y$ be the cokernel of he map
$X_d \to X_{d-1}$ and let $\la \in\wedge^+ (n,r)$ be such that
$K_\la V$ is a factor of $X_d$ but not of $X_{d-1}$. From the exact
sequence $0\to X_d \to X_{d-1}\to Y\to 0$ we obtain the exact
sequence
\[
    {\rm Hom} (X_{d-1},L_{\tilde{\la}}V)\lto {\rm Hom}(X_d
    ,L_{\tilde{\la}}V)\lto {\rm Ext}^1 (Y,L_{\tilde{\la}}V)\lto {\rm
    Ext}^1 (X_{d-1} ,L_{\tilde{\la}}V).
\]
By Proposition 2.1 ii) we have ${\rm Ext}^1
(X_{d-1},L_{\tilde{\la}}V)=0$ and by the assumption on $\la$, ${\rm
Hom}(X_{d-1},L_{\tilde{\la}}V)=0$. Thus
\[
    {\rm Ext}^1 (Y,L_{\tilde{\la}}V)={\rm Hom}(X_d ,L_{\tilde{\la}}
    V)\neq 0.
\]
By dimension shifting on $0\to X_d \to \cdots \to X_0 \to X\to 0$ we
have ${\rm Ext}^d (X, L_{\tilde{\la}}V)={\rm Ext}^1
(Y,L_{\tilde{\la}}V)$ so that ${\rm Ext}^d (X, L_{\tilde{\la}}
V)\neq 0$. Hence $wfd (X)\ge d$ by Proposition 2.1 ii).\\
\end{proof}

From the previous Lemma and the complexes $\K_* (r,1,V)$ of Theorem
1.4 we obtain a lower bound for $wfd(S(n,r))$ if $p=2$ and $r$ is
even valid for any $n$.

\begin{corollary}%\label{co2.3}
Suppose $p=2$ and $r$ is even, $r=2r'$. Then
\[
    wfd((D_{r'}V)^{(1)})=\frac{r'}{2} \ \text{ and } \ wfd(S(n,r))\ge
    \frac{r'}{2}.
\]
\end{corollary}

\noindent\textbf{2.3} We consider here the Schur algebra $S(2,r)$
and give an alternative proof of the following result due to Parker
\cite{P1}.

\begin{theorem}%\label{th2.4}
1) If $p=2$, then
\[
    wfd (S(2,r))=\left\{\begin{array}{ll}
    \frac{r}{2},&\text{if $r$ is even}\\
    \lfloor \frac{r}{4}\rfloor, &\text{if $r$ is odd.}
    \end{array}\right.
\]
2) If $p\ge 3$, then
\[
    wfd (S(2,r))=\left\lfloor \frac{r}{p}\right\rfloor .
    \ \ \hspace{2cm}
\]
\end{theorem}

In order to prove this result, we will establish bounds on $wfd
(S(2,r))$ by exhibiting suitable Weyl resolutions.

Suppose $\la =(a,b)\in \wedge^+ (2)$ and let $c=c(\la)=a-b$. First
we observe that since $K_{(b,b)}V$ is a $1$-dimensional irreducible
representation of highest weight $(b,b)$, we have $K_\la
V=K_{(b,b)}V\otimes D_c V$ and for each $m\ge 1$, $K_\la
V^{(m)}=K_{(p^m b,p^m b)}V\otimes D_c V^{(m)}$.

\bigskip\noindent\textbf{2.3.1 Upper bound for any $p$}

\begin{lemma}%\label{le2.5}
Let $\la =(a,b)\in \wedge^+ (2)$ and $c=a-b$. Then for all $m\ge 1$,
we have $wfd(K_\la V^{(m)})\le p^{m-1}c$.
\end{lemma}

\begin{proof}
We use induction on $m$. For $m=1$, it follows from the above
observation and Theorem 1.4 that the complex $C_*
=K_{(pb,pb)}V\otimes \K_* (pc,1,V)$ is a Weyl resolution of its
homology which is $K_\la V^{(1)}$. Since the length of $C_*$ is $c$,
we have $wfd(K_\la V^{(1)})\le c$.

Suppose $m\ge 2$ and let $C_{*}^{(m)}$ be the $m$-th Frobenius twist
of the complex $C_*$. Then $C_{*}^{(m)}$ has length $c$ and is a
resolution of $K_\la V^{(m+1)}$. From the definition of $\K_*
(pc,1,V)$, it follows that $C_{i}^{(m)} =K_{\la_i}V^{(m)}$ for some
partition $\la_i$ satisfying $c(\la_i )\le pc-2i$. By the inductive
hypothesis, $wfd (C_{i}^{(m)})\le p^{m-1}(pc -2i)$. From Proposition
2.1 vi), $wfd(K_\la V^{(m+1)})\le \max \{ p^{m-1} (pc-2i)+i|
i=0,\cdots, c\} =p^m c$.
\end{proof}

Now if $\la =(a,b)\in \wedge^+ (2,r)$, $c=a-b$ and $c=c_0
p^{n_0}+\cdots + c_s p^{n_s}$ is the base $p$ expansion of $c$, then
by Steinberg's tensor product theorem
\begin{align*}
    L(\la) &= L(b,b)\otimes L(c_0 )^{(n_0 )}\otimes \cdots \otimes
    L(c_s )^{(n_s )}\\
    &=K_{(b,b)}V\otimes D_{c_0}V^{(n_0 )}\otimes \cdots \otimes
    D_{c_s} V^{(n_s )}.
\end{align*}
By Proposition 2.1 iv) and the above lemma we have
\[
    wfd (L(\la ))\le p^{n_0 -1} c_0 +\cdots + p^{n_s -1}c_s =
    \frac{c}{p} \le \left\lfloor \frac{r}{p}\right\rfloor .
\]
By Proposition 2.1 v), we have $wfd (S(2,r))\le \lfloor
\frac{r}{p}\rfloor$.

\bigskip\noindent\textbf{2.3.2 Lower bound for $p\ge 3$}

\medskip Let $r=pr_1 +r_0$, where $0\le r_0 <p$, and $r_0 =2s+t$,
where $0\le t<2$. By Theorem 1.4, the complex $K_{(s,s)} V\otimes
\K_* (pr_1 +t,1+t,V)$ is a Weyl resolution of its homology. This
complex has length $r_1 =\lfloor\frac{r}{p}\rfloor$ and satisfies
the hypothesis of Lemma 2.2 so that its homology has Weyl filtration
dimension $\lfloor\frac{r}{p}\rfloor$. Hence $wfd (S(2,r) )\ge
\lfloor\frac{r}{p}\rfloor$.

From (2.3.1) and (2.3.2) we have $wfd
(S(2,r))=\lfloor\frac{r}{p}\rfloor$ if $p\ge 3$.

\bigskip\noindent\textbf{2.3.3 Better upper bound for $p=2$ and $r$ odd}

\medskip First we observe that if $X$ has a filtration with factors
various $K_\la V^{(1)}$, then $V\otimes X$ has a Weyl filtration.
This follows from the fact that $V\otimes K_\la V^{(1)}$ is a Weyl
module (in the special case under consideration, $p=n=2$, the
natural module $V$ is a Steinberg module).

\begin{lemma}%\label{le2.6}
Suppose $\la_i =(a_i ,b_i )\in \wedge^+ (2)$, $i=1,\dots, s$ and let
$c_i =a_i -b_i$. Then for all positive integers $n_1 ,\dots, n_s$ we
have
\[
    wfd (V\otimes K_\la V^{(n_1 )} \otimes \cdots \otimes
    K_{\la_s}^{(n_s )}V)\le c_1 2^{n_1 -2} +\cdots + c_s 2^{n_s -2}.
\]
\end{lemma}

\begin{proof}
We use induction on $\sum\limits^{s}_{i=1}n_i$. If $n_1 =\cdots =n_s
=1$, then by Proposition 2.1 iii) $K_{\la_1}V^{(1)}\otimes \cdots
\otimes K_{\la_s}V^{(1)}$ has a filtration with quotients various
$K_\la V^{(1)}$. By the observation mentioned before the Lemma, it
follows that $V\otimes K_{\la_1}V^{(1)}\otimes \cdots \otimes
K_{\la_s}V^{(1)}$ has a Weyl filtration and thus its Weyl filtration
dimension is 0.

The complex $C_* =K_{(2b_1 ,2b_1 )}V\otimes \K_* (2c_1 ,1,V)$ is a
Weyl resolution of $K_{\la_1}V^{(1)}$ of length $c_1$. Hence the
complex $\bar{C}_* =V\otimes C_{*}^{(n_1 )}\otimes K_{\la_2}V^{(n_2
)} \otimes \cdots \otimes K_{\la_s}V^{(n_s )}$ is a resolution of
$V\otimes K_{\la_1}V^{(n_1 +1)} \otimes K_{\la_2}V^{(n_2 )}\otimes
\cdots \otimes K_{\la_s}V^{(n_s )}$ of length $c_1$.

From the definition of $\K_* (2c_1 ,1,V)$, it follows that
$C_{i}^{(n_1 )}=K_{\mu_i}V^{(n_1 )}$ where $\mu_i \in \wedge^+ (2)$
satisfies $c(\mu_i )\le 2c_1 -2i$. By applying Proposition 2.1 vi)
to $\bar{C}_*$ and using the inductive hypothesis, we have
\begin{align*}
    wfd (V &\otimes K_{\la_1}V^{(n_1 +1)}\otimes K_{\la_2}V^{(n_2
    )}\otimes\cdots\otimes K_{\la_s}V^{(n_s )})\le\\
    &\le \max \left\{ (2c_1 -2i)2^{n_1 -2}+c_2 2^{n_2 -2} +\cdots +
    c_s 2^{n_s -2}+i | i=0,\dots ,c_1 \right\} \\
    &=c_1 2^{n_1 -1}+c_2 2^{n_2 -2} +\cdots + c_2 2^{n_s -2}.
\end{align*}
\end{proof}
Now suppose $r$ is odd, $\la =(a,b)\in \wedge^+ (2,r)$, $c=b-a$ and
$c=1+2^{n_1} +\cdots + 2^{n_s}$ is the base 2 expansion of $c$. Then
$L(\la)=K_{(b,b)}V\otimes V\otimes V^{(n_1 )}\otimes \cdots \otimes
V^{(n_s )}$ and by Proposition 2.1 iv) and the above Lemma,
$wfd(L(\la)) \le wfd(V\otimes V^{(n_1 )}\otimes \cdots \otimes
V^{(n_s )})\le 2^{n_1 -2}+\cdots +2^{n_s -2}=\frac{c-1}{4}\le
\lfloor\frac{r}{4}\rfloor$. Hence $wfd S(2,r)\le
\lfloor\frac{r}{4}\rfloor$.

If $r$ is even, we have he upper bound $wfd (S(2,r))\le \frac{r}{2}$
from (2.3.1).

\bigskip\noindent\textbf{2.3.4 Lower bound for $p=2$}

If $r$ is even, we have $wfd (S(2,r))\ge \frac{r}{2}$ by Corollary
2.3.

Suppose $r$ is odd. If $r\equiv 1\md 4$, then from the complex
$\mathbb{M}_* (r,2,V)$ (cf. subsection 1.5) and Lemma 2.1 we see
that the homology of $\mathbb{M}_* (r,2,V)$ has Weyl filtration
dimension $\lfloor\frac{r}{4}\rfloor$. Hence $wfd (S(2,r))\ge
\lfloor\frac{r}{4}\rfloor$. If $r\equiv 3\md 4$ we obtain the same
conclusion by considering the complex $\N_{*}(r,2,V)$ of subsection
1.5.

From (2.3.3) and (2.3.4), we obtain the value of $wfd (S(2,r))$ in
Theorem 2.3.

%D_{a+b-\ell}V\otimes D_{\ell}V
%%%%%%%%%%%%%%%%%%%%%%%%%%%%%%%%%%%%%%%%%%%%%%%%%%%%%%%%%%%%%%%%%%

\medskip\noindent
Mihalis Maliakas\\
Department of Mathematics\\
University of Athens\\
Panepistimiopolis 15784\\
Greece\\
email: mmaliak@math.uoa.gr


\begin{thebibliography}{99}
\bibitem{A-B-W}%1
K. Akin, D. Buchsbaum and J. Weyman, \textit{Schur functors and
Schur complexes}, Adv. in Math. 44(1982), 207--278.

\bibitem{C-P}%2
R.W. Carter and M.T.J. Payne, \textit{On homomorphisms between Weyl
modules and Specht modules}, Math. Proc. Camb. Phil. Soc. 87(1980),
419--425.\\[-25pt]

\bibitem{D}%3
S. Donkin, \textit{Finite resolutions of modules for reductive
algebraic groups}, J. Algebra 101(1986), 473--488.\\[-25pt]

\bibitem{F-P}%4
E.M. Friedlander and B.J. Parshall, \textit{Cohomology of Lie
algebras and algebraic groups}, Amer. J. Math. 108(1986), 235--253.\\[-25pt]

\bibitem{G}%5
J.A. Green, \textit{Polynomial Representations of $GL_n$}, $2^{\rm
nd}$ Edition, LNM 830, Springer-Verlag, Berlin Heiderberg,
2007.\\[-25pt]

\bibitem{J}%6
J.C. Jantzen, \textit{Representations of Algebraic Groups}, $2^{\rm
nd}$ Edition, Math. Surveys and Monographs 107, AMS, 2003.\\[-25pt]

\bibitem{M}%7
I.G. Macdonald, \textit{Symmetric Functions and Hall Polynomials},
$2^{\rm nd}$ Edition, Clarendon Press, Oxford, 1995.\\[-25pt]

\bibitem{P1}%8
A.E. Parker, \textit{The global dimension of Schur algebras for
$GL_2$ and $GL_3$}, J. Algebra 241(2001), 340--378.\\[-25pt]

\bibitem{P2}%9
A.E. Parker, \textit{The global filtration dimension of Weyl modules
for a linear algebraic group}, J. reine angew. Math. 562(2003),
5--21.\\[-25pt]

\bibitem{W}%10
J. Wang, \textit{Sheaf cohomology on $G/B$ and tensor product of
Weyl modules}, J. Algebra 77(1982), 162--185.\\[-25pt]
\end{thebibliography}
\end{document}